\documentclass{amsart}
\usepackage{amsmath}
\usepackage{amssymb}
\usepackage{amsthm}
\usepackage{mathrsfs}

\usepackage{graphicx}

\newtheorem{example}{Example}[section]
\newtheorem{theorem}{Theorem}[section]
\newtheorem{corollary}{Corollary}[section]
\newtheorem{lemma}{Lemma}[section]

\newcommand{\RR} {\mathbb R}

\newcommand{\NN} {\mathbb N}

\newcommand{\RRR} {\mathcal R}

\newcommand{\MMM} {\mathcal M}

\newcommand{\ra} {\rightarrow}
\theoremstyle{definition}

\theoremstyle{remark}

\numberwithin{equation}{section}

\begin{document}

\title{Real compactness via real maximal ideals of $B_1(X)$}

\author{A. Deb Ray}
\address{Department of Pure Mathematics, University of Calcutta, 35, Ballygunge Circular Road, Kolkata - 700019, INDIA} 
\email{adrpm@caluniv.ac.in}

\author{Atanu Mondal}
\address{Department of Commerce (E), St. Xavier's college, 30, Mother Teresa sarani, Kolkata - 700016, INDIA}
\email{atanu@sxccal.edu}

\begin{abstract}
In this paper, constructing a class of ideals of $B_1(X)$ from proper ideals of $C(X)$ a one-one correspondence between the class of real maximal ideals of $C(X)$ and those of $B_1(X)$ is established.  The collection of all real maximal ideals of $B_1(X)$ with hull-kernel topology is proved to be homeomorphic to the space of real maximal ideals of $C(X)$ endowed with a topology finer than the subspace topology induced from its structure space. It is also proved that a Tychonoff space is real compact if and only if every real maximal ideal of $B_1(X)$ is fixed. As a consequence, within the class of real compact $T_4$ spaces whose points are $G_\delta$, $B_1(X) = {B_1}^*(X)$ if and only if $X$ is finite. 
\end{abstract}

\maketitle

\section{Introduction}
\noindent  The pointwise limit function of a sequence $\{f_n\}$ of real valued continuous functions defined on a topological space $X$ is well known as a Baire class-one (or, Baire one) function. Several characterizations of Baire one functions defined on metric spaces were obtained by different mathematicians \cite{JPFEACAPR}, \cite{DZ}. Inspired by the research on rings of continuous functions, in \cite{AA1}, we introduced two rings $B_1(X)$ and ${B_1}^*(X)$ consisting respectively of real valued Baire one functions and bounded Baire one functions on a topological space $X$.  It has been observed in \cite{AA1} that $B_1(X)$ is a commutative lattice ordered ring with unity with respect to the usual addition and multiplication of functions and it is an over-ring of the ring $C(X)$ of all real valued continuous functions on $X$. The main purpose of this paper is to establish that the class of all real maximal ideals of $B_1(X)$ is in one-one correspondence with the class of all  real maximal ideals of $C(X)$. In fact, defining a sort of `extension' of an ideal $I$ of  $C(X)$ in $B_1(X)$ (denoted by $I_B)$, we show that the contraction $I_B\cap C(X)$ coincides with $I$ if and only if $I$ is a real maximal ideal in $C(X)$. That $I_B$ is a real maximal ideal in $B_1(X)$, for any real maximal ideal $I$ in $C(X)$ determines the said one-one correspondence. Moreover, it is proved that the collection $\mathcal{RM}(B_1(X))$ of all real maximal ideals of $B_1(X)$ with the subspace topology (i.e., the  hull-kernel topology) of the structure space of $B_1(X)$ is homeomorphic to the collection $\mathcal{RM}(C(X))$ of all real maximal ideals of $C(X)$ equipped with a topology finer than the subspace topology of the structure space of $C(X)$.\\\\
\noindent The class of topological spaces on which every Baire one function is bounded, is yet to be determined completely. In Section 3, we prove a necessary and sufficient condition for $B_1(X)$ to coincide with ${B_1}^*(X)$ within the class of all $T_4$, real compact spaces whose points are $G_\delta$.\\\\
\noindent To make this paper self-sufficient, we now introduce some known terminologies and facts. A zero set of $f \in B_1(X)$ is defined as usual by a set of the form $Z(f) = \{x\in X \ : \ f(x) = 0\}$. As an analogue of $z$-filter (or $z$-ultrafilter) on $X$, we introduced in \cite{AA2} the $Z_B$-filter (or respectively, $Z_B$-ultrafilter) on $X$, thereby investigating the duality between ideals (maximal ideals) in $B_1(X)$ and $Z_B$-filters (respectively, $Z_B$-ultrafilters) on $X$. The above mentioned duality existing between ideals in $B_1(X)$ and $Z_B$-filters on $X$ is manifested by the fact that if $I$ is an ideal in $B_1(X)$ then $Z_B[I] = \{Z(f) \ : \ f\in B_1(X)\}$ is a $Z_B$-filter on $X$ and dually for a $Z_B$-filter $\mathscr F$ on $X$, $Z_B^{-1}[\mathscr F] = \{f\in B_1(X) \ : \ Z(f) \in \mathscr F\}$ is a proper ideal in $B_1(X)$. The assignment $M \mapsto Z_B[M]$ is a bijection from the set of all maximal ideals in $B_1(X)$ and to the family of all $Z_B$-ultrafilters on $X$. In the same paper \cite{AA2}, defining suitably the residue class fields $B_1(X) / M$, the concept of real and hyper-real maximal ideals of $B_1(X)$ was introduced. A maximal ideal $M$ of $B_1(X)$ is called real \cite{AA2} if $B_1(X)/M \cong \RR$ and in such case $B_1(X)/M$ is called real residue class field. If $M$ is not real then it is called hyper-real \cite{AA2} and $B_1(X)/M$ is called hyper-real residue class field. Considering the structure space $\MMM(B_1(X))$ of $B_1(X)$, i.e., the collection $\MMM(B_1(X))$ of all maximal ideals of $B_1(X)$ with respect to the hull-kernel topology, we get the subspace topology on the collection $\RRR\MMM(B_1(X))$ of all real maximal ideals of $B_1(X)$. We show that the subspace topology via the aforesaid bijection induces a topology on $\RRR\MMM(C(X))$ of all real maximal ideals of $C(X)$ which is finer than the hull-kernel topology on $\RRR\MMM(C(X))$. \\\\
\noindent Throughout this paper, we consider $X$ as any Hausdorff topological space unless stated otherwise. A function $f : X \ra \RR$ is called a Baire one function if there exists a sequence of continuous functions $\{f_n\}$ such that for all $x\in X$, $\{f_n(x)\}$ converges to $f(x)$. We use the notation $f_n \xrightarrow []{p. w.} f$ to denote  $\{f_n\}$ pointwise converges to $f$. 

\begin{theorem}\label{thm1.1}\textnormal{[Theorem 5.14 of \cite{GJ}]}
	For a maximal ideal $M$ of $C(X)$ the following statements are equivalent:
	\begin{enumerate}
		\item $M$ is a real maximal ideal.
		\item $Z[M]$ is closed under countable intersection.
		\item $Z[M]$ has countable intersection property.
	\end{enumerate}
\end{theorem}
\noindent An analogue of this theorem in the context of Baire one functions is the following:
\begin{theorem}\label{RMaxB}\textnormal{[Theorem 4.26 of \cite{AA2}]}
	For a maximal ideal $M$ of $B_1(X)$ the following statements are equivalent:
	\begin{enumerate}
		\item $M$ is a real maximal ideal.
		\item $Z_B[M]$ is closed under countable intersection.
		\item $Z_B[M]$ has countable intersection property.
	\end{enumerate}
\end{theorem}

\section{A one-one correspondence between the real maximal ideals of $C(X)$ and the real maximal ideals of $B_1(X)$}

\noindent It is easy to observe that for each proper ideal $I$ of $C(X)$, $I_B = \{f\in B_1(X) \ : \ \exists \{f_n\} \subseteq I \textnormal{ such that } f_n \xrightarrow []{p. w.} f\}$ is an ideal of $B_1(X)$ such that $I \subseteq I_B \cap C(X)$.
\noindent As a natural example, we obtain that $(M_p)_B$ is a fixed maximal ideal of $B_1(X)$ \cite{AA2}, where $M_p=\{f \in C(X): f(p)=0\}$ is a fixed maximal ideal of $C(X)$ \cite{GJ} and this example prompts us to prove a more general result later in this section.
\begin{example}\label{FixedMax1}
For each $p\in X$, ${(M_p)}_B = \widehat{M}_p \equiv \{f \in B_1(X): f(p)=0\}$.\\
For each $f \in {(M_p)}_B$, there exists $\{f_n\} \subseteq M_p$ such that $f_n \xrightarrow []{p. w.} f$. This implies $f_n(p) = 0$, for all $n \in \NN$ and therefore, $f(p) = 0$. i.e., ${(M_p)}_B \subseteq \widehat{M}_p$. On the other hand, $f \in \widehat{M}_p$ implies $f(p) = 0$. Since $f \in B_1(X)$, there exists a sequence $\{g_n\} \subseteq C(X)$ such that $g_n \xrightarrow []{p. w.} f$. Define$f_n = g_n - g_n(p)$, for all $n \in \NN$. Clearly, $f_n(p) = 0$, for all $n\in \NN$. Also, $f_n \xrightarrow []{p. w.} f$. Therefore, $f\in {(M_p)}_B$ and $\widehat{M}_p \subseteq {(M_p)}_B$. Hence, ${(M_p)}_B = \widehat{M}_p$.
\end{example}

\begin{theorem}
If $I$ is an absolutely convex ideal in $C(X)$ then $I_B$ is an absolutely convex ideal in $B_1(X)$.
\end{theorem}
\begin{proof}
We first prove that $I_B$ is a convex ideal in $B_1(X)$. If so, then $f \in I_B$ implies that there is $\{f_n\} \subseteq I$ such that $f_n \xrightarrow []{p. w.} f$ and hence, $|f_n| \xrightarrow []{p. w.} |f|$. As $I$ is absolutely convex, $\{|f_n|\} \subseteq I$ which ensures $|f| \in I_B$. In such case, $I_B$ becomes absolutely convex.\\
Let $f, g \in B_1(X)$ such that $0 \leq f \leq g$ and $g \in I_B$. Then there is a sequence $\{f_n\}$ in $C(X)$ and $\{g_n\} \subseteq I$ such that $f_n \xrightarrow []{p. w.} f$ and $g_n \xrightarrow []{p. w.} g$. Choosing $h_n = f_n \wedge g_n$, we observe the following:
\begin{enumerate}
	\item[(i)] $h_n \xrightarrow []{p. w.} f \wedge g = f$.
	\item[(ii)] For each $n \in \mathbb{N}$, $0 \leq h_n \leq g_n$ and $g_n \in I$ implies that $h_n \in I$ (since, $I$ is absolutely convex).
\end{enumerate}
Hence, $f\in I_B$. 
\end{proof}
\noindent For any proper ideal $I$ of $C(X)$, it is clear that $I \subseteq I_B \cap C(X)$. In the following theorem we show that the equality holds precisely for the class of all real maximal ideals of $C(X)$.

\begin{theorem}\label{RealIdeal1}
$M \in \RRR\MMM(C(X))$ if and only if $M = M_B \cap C(X)$.
\end{theorem}
\begin{proof}
Let $M$ be a real maximal ideal of $C(X)$. Clearly, $M \subseteq M_B \cap C(X)$. Now let $g\in M_B \cap C(X)$. There exists $\{g_n\} \subseteq M$ such that $g_n \xrightarrow []{p. w.} g$. Since $M$ is real and $g_n \in M$, for all $n\in \NN$, $\bigcap\limits_{n=1}^{\infty} Z(g_n) \in Z[M]$. Also, $\bigcap\limits_{n=1}^{\infty} Z(g_n) \subseteq Z(g)$. Hence, $Z(g) \in Z[M]$. By maximality of $M$ it follows that $g\in M$. Therefore, $M_B \cap C(X) \subseteq M$ and it implies that $M = M_B \cap C(X)$.  \\
Conversely, let $M$ be a maximal ideal of $C(X)$ such that $M = M_B \cap C(X)$.\\
Consider any countable family $\{Z(g_n) : n \in \NN\}$ of $Z[M]$. By maximality of $M$, $g_n \in M$, for all $n \in \NN$. \\
We now construct a sequence $\{s_n\}$ as follows : $s_n = \sum\limits_{j = 1}^n \left( \frac{1}{2^j} \wedge |g_j|\right)$, for each $n\in \NN$. Certainly, for each $j$, $Z(g_j) = Z\left(\frac{1}{2^j} \wedge |g_j|\right)$ implies that $\frac{1}{2^j} \wedge |g_j| \in M$. $M$ being an ideal, finite sum of each such member will also lie within $M$ and hence, $s_n \in M$, for all $n\in \NN$.\\
$s = \sum\limits_{n=1}^\infty \left(\frac{1}{2^n} \wedge |g_n|\right)$ is the uniform limit of the sequence $\{s_n\}$ of continuous functions and therefore, $s \in C(X)$. Again $\{s_n\} \subseteq M$ ensures that $s \in M_B \cap C(X) = M$. So, $Z(s) \neq \emptyset$. Following the arguments  used in {\rm (1.14 (a) of \cite {GJ})} we obtain $\bigcap\limits_{n=1}^{\infty} Z(g_n) = Z(s) \neq \emptyset$. Hence, by Theorem \ref{thm1.1} $M$ is real.
\end{proof}
\noindent That $M_B$ is not even a proper ideal of $B_1(X)$ when $M$ is hyper-real in $C(X)$ is observed in the next theorem.
\begin{theorem}
If $M$ is a hyper-real maximal ideal in $C(X)$ then $M_B = B_1(X)$.
\end{theorem}
\begin{proof}
If $M$ is hyper-real then by Theorem~\ref{RealIdeal1} $M_B \cap C(X) \neq M$. But for any ideal $I$, $I \subseteq I_B \cap C(X)$ holds. So, $M \subsetneq M_B \cap C(X)$. Since $M$ is maximal, $M_B \cap C(X) = C(X)$. Hence, $C(X) \subseteq M_B$, i.e., $1 \in M_B$. Therefore, $M_B = B_1(X)$. 
\end{proof}
\begin{theorem}\label{RealMax2}
If $M\in \RRR\MMM(C(X))$ then $M_B \in \RRR\MMM(B_1(X))$. 
\end{theorem} 
\begin{proof}
Let $f\in B_1(X) \setminus M_B$. Consider the ideal J generated by $M_B \cup \{f\}$.\\
$f \in B_1(X)$ implies that there exists $\{f_n\} \subseteq C(X)$ such that $f_n \xrightarrow []{p. w.} f$.\\
Since $M$ is a real maximal ideal in $C(X)$, for each $f_n\in C(X)$, there exists some $r_n \in \RR$ such that $M(f_n) = M(r_n)$ and so, $f_n = r_n$ on $Z_n = Z(f_n-r_n) \in Z[M]$. Since $Z = \bigcap\limits_{n=1}^{\infty} Z_n \in Z[M]$, each $f_n= r_n$ on $Z$. As a consequence, $f$ is constant (say, r) on $Z \in Z[M] \subseteq Z_B[M_B]$, where $r = \lim\limits_{n \rightarrow \infty} r_n$. \\
Since, $Z \subseteq Z(f_n-r_n)$ implies that $Z(f_n-r_n)\in Z[M]$ and $M$ is a $z$-ideal in $C(X)$, $f_n-r_n \in M$. By definition of $M_B$, $f-r \in M_B$. Since $f\notin M_B$, $r\neq 0$. But, $r =f - (f-r) \in J$ and $r\neq 0$ implies that $J = B_1(X)$. Hence, $M_B$ is a maximal ideal of $B_1(X)$ such that $f-r \in M_B$. i.e., $M_B(f) = M_B(r)$, for some $r\in \RR$. If $f \in M_B$ then $M_B(f) = M_B(0)$ and this proves that $M_B \in \RRR\MMM(B_1(X))$.
\end{proof}

\begin{theorem}\label{RealMax3}
If $\widehat{M} \in \RRR\MMM(B_1(X))$ then $\widehat{M}\cap C(X) \in \RRR\MMM(C(X))$.
\end{theorem}
\begin{proof}
Let $\widehat{M} \in \RRR\MMM(B_1(X))$. For each $f\in B_1(X)$, there exists $r_f \in \RR$ such that $f - r_f \in \widehat{M}$. In particular, for any $f \in C(X)$, there is $r_f \in \RR$ such that $f - r_f \in \widehat{M}$. So, $f - r_f \in \widehat{M}\cap C(X) = M$ (say). Define a function $\phi : C(X)/M \ra \RR$ by $M(f) \mapsto r_f$, whenever $f - r_f \in M$. We claim that $\phi$ is an isomorphism.\\
$M(f) = M(g) \Leftrightarrow f-g \in M$. If $\phi(M(f)) = r_f$ and $\phi(M(g)) = r_g$ then $f-r_f, g-r_g \in M$. i.e., $(f - g) - (r_f-r_g) \in M$. Since, $f-g \in M$ and $M$ is an ideal, it follows that $r_f - r_g \in M$ - a contradiction to the fact that $M$ is proper, unless $r_f - r_g = 0$. Hence $\phi$ is well defined. \\
$\phi(M(f)) = \phi(M(g))$ implies that $r_f = r_g$ where $f-r_f, g - r_g \in M$. Therefore, $f - g = (f - r_f) - (g - r_g) \in M$ which in turn gives $M(f) = M(g)$, proving $\phi$ to be one-one. \\
The function $\phi$ is clearly onto, as $\phi(M(r)) = r$, for each $r \in \RR$. By routine arguments we easily see that $\phi$ is indeed a ring homomorphism. Hence, $\phi$ is a ring isomorphism and therefore, $M \in \RRR\MMM(C(X))$.
\end{proof}
\begin{corollary}\label{Corollary1}
If $\widehat{M} \in \RRR\MMM(B_1(X))$ then $(\widehat{M} \cap C(X))_B = \widehat{M}$.
\end{corollary}
\begin{proof}
As $\widehat{M} \in \RRR\MMM(B_1(X))$, $\widehat{M} \cap C(X) \in \RRR\MMM(C(X))$ (by Theorem~\ref{RealMax3}). Using Theorem~\ref{RealMax2}, $(\widehat{M}\cap C(X))_B \in \RRR\MMM(B_1(X))$. Since $(\widehat{M}\cap C(X))_B$ is a maximal ideal, it is enough to show that $(\widehat{M}\cap C(X))_B \subseteq \widehat{M}$.\\
Let $g\in (\widehat{M}\cap C(X))_B$ . Then there exists $\{g_n\} \subseteq \widehat{M}\cap C(X)$ such that $g_n \xrightarrow []{p. w.} g$. So, $Z(g) \supseteq \bigcap\limits_{i=1}^\infty Z(g_n)$. As $Z_B[\widehat{M}]$ is a $Z_B$-ultrafilter and $\widehat{M}$ is real, it follows that $Z(g) \in Z_B[\widehat{M}]$. Hence, $g\in \widehat{M}$ and therefore $(\widehat{M}\cap C(X))_B \subseteq \widehat{M}$. 
\end{proof}

\noindent In view of Corollary~\ref{Corollary1}, Theorem~\ref{RealMax2} and Theorem~\ref{RealMax3}, we get a one-one correspondence between $\RRR\MMM(C(X))$ and $\RRR\MMM(B_1(X))$
\begin{theorem}\label{bij}
	If $\psi : \RRR\MMM(C(X)) \rightarrow \RRR\MMM(B_1(X))$ is defined by $M \mapsto M_B$ then $\psi$ is a bijection.
\end{theorem}
\begin{proof}
	Let $\widehat{M}$ be any member of $\RRR\MMM(B_1(X))$. Therefore, by Corollary \ref{Corollary1} we get $(\widehat{M} \cap C(X))_B=\widehat{M}$, where $\widehat{M} \cap C(X) \in \RRR\MMM(C(X))$ (By Theorem \ref{RealMax3}). Hence, for $\widehat{M} \in \RRR\MMM(B_1(X)),$ we get $\widehat{M} \cap C(X) \in \RRR\MMM(C(X))$ such that $\psi(\widehat{M} \cap C(X))=\widehat{M}.$ This proves that $\psi$ is surjective.\\
	To show that, $\psi$ is injective we assume $\psi(\widehat{M})=\psi(\widehat{N})$. This implies $(\widehat{M})_B=(\widehat{N})_B$. Now by applying Theorem \ref{RealIdeal1}, we get $\widehat{M}= (\widehat{M})_B \cap C(X)= (\widehat{N})_B \cap C(X)= \widehat{N}$. Therefore, $\psi$ is surjective and hence, it is a bijection.
\end{proof}

\begin{corollary}\label{Real}
$|\RRR\MMM(C(X))| = |\RRR\MMM(B_1(X))|$.
\end{corollary} 
\noindent It is well known that $\{ \widehat{{\mathscr M}}_f \ : \ f\in B_1(X)\}$, where each $\widehat{{\mathscr M}}_f  = \{M \in \MMM(B_1(X)) : f\in M\}$, forms a base for closed sets for the hull-kernel topology on $\MMM(B_1(X))$ and certainly $\RRR\MMM (B_1(X))$ is a subspace of $\MMM(B_1(X))$. In the following theorem we show that the bijection $\psi$ obtained above  becomes a homeomorphism if  $\RRR\MMM(C(X))$ is endowed with a finer topology than the subspace topology induced from the hull-kernel topology of $\MMM(C(X))$.

\begin{theorem}
Let $(X, \tau)$ be a Tychonoff space. Then for each $f\in B_1(X)$, the collection ${{\mathscr M}_f}^*= \{M \in \RRR\MMM(C(X)) \ : \ f\in M_B\}$ forms a base for closed sets for some topology $\sigma$ on $\RRR\MMM(C(X))$ which is finer than the subspace topology of the structure space of $\MMM(C(X))$. Moreover, $\psi : (\RRR\MMM(C(X)), \sigma) \ra \RRR\MMM(B_1(X))$ given by $M \mapsto M_B$ is a homeomorphism.
\end{theorem}
\begin{proof}
To prove $\mathscr{B}^*=\{{{\mathscr M}_f}^*: f \in B_1(X)\}$ forms a base for closed sets for some topology $\sigma$ on $\RRR\MMM(C(X))$, it is enough to show that $\emptyset \in \mathscr{B}^*$ and $\mathscr{B}^*$ is closed under finite union. It is easy to observe that, $\emptyset = {{\mathscr M}_1}^* \in \mathscr{B}^*$. Now let ${{\mathscr M}_f}^*$, ${{\mathscr M}_g}^* \in $ $\mathscr{B}^*$, for some $f,g \in B_1(X)$. Take any $M \in {{\mathscr M}_f}^* \cup {{\mathscr M}_g}^*$. Therefore, $gf \in M_B$ and $M \in {{\mathscr M}_{fg}}^*$. This implies ${{\mathscr M}_f}^* \cup {{\mathscr M}_g}^* \subseteq {{\mathscr M}_{fg}}^*$. On the other hand, if we take any memeber $M \in {{\mathscr M}_{fg}}^*$ then we get $fg \in M_B$. Now by Theorem \ref{RealMax2} $M_B$ is a maximal ideal and hence $f \in M_B$ or $g \in M_B$, i.e., $M \in {{\mathscr M}_{f}}^* \cup {{\mathscr M}_{g}}^*$. So ${{\mathscr M}_{fg}}^* \subseteq {{\mathscr M}_f}^* \cup {{\mathscr M}_g}^*$. This proves that ${{\mathscr M}_f}^* \cup {{\mathscr M}_g}^*= {{\mathscr M}_{fg}}^*$ and hence, $\mathscr{B}^*$ is closed under finite union.\\
Now to prove that $\psi$ is a homeomorphism, we need to show $\psi$ is bijective and exchanges the basic closed sets of $(\RRR\MMM(C(X)), \sigma)$ and $\RRR\MMM(B_1(X))$. The map $\psi$ is bijective is already proved in Thoeorem \ref{bij}. Now for any $f \in B_1(X),$ $\psi({{\mathscr M}_f}^*)=\{\psi(M): f \in M_B\}=\{M_B:f \in M_B\}=\{N \in \RRR\MMM(B_1(X)):f \in N\}=\widehat{{\mathscr M}}_f \bigcap \RRR\MMM(B_1(X))$, which is a basic closed set of $\RRR\MMM(B_1(X))$ for the subspace topology induced from the hull-kernel topology on $\MMM(B_1(X))$.
As $\psi$ exchanges the basic closed sets, it is a homeomorphism.
\end{proof}

\noindent Before we conclude this section, we show that an injective map exists from $\mathcal{H}(C(X))$ into $\mathcal{H}(B_1(X))$, where $\mathcal{H}(C(X))$ and $\mathcal{H}(B_1(X))$ represent the collections of all hyper-real maximal ideals in $C(X)$ and $B_1(X)$ respectively. In what follows, we use the notation $I^*$ for the ideal of $B_1(X)$ generated by the subset $I$ of $B_1(X)$ and $m(I^*)$ for its maximal extension. The next theorem ensures that the ideal of $B_1(X)$ generated by a proper ideal of $C(X)$ is indeed proper, so that it has a maximal extension, say $m(I^*)$.
\begin{theorem}
For any proper ideal $I$ of $C(X)$, $I^*$ is a proper ideal of $B_1(X)$, where $I^*$ denotes the ideal of $B_1(X)$ generated by $I$ as a subset of $B_1(X)$.
\end{theorem}
\begin{proof}
If possible let, $I^*$ be not proper. Then $I^* = B_1(X)$ and hence $\textbf{1}$ (the constant function with value 1) can be written as $\mathbf{1}= \sum\limits_{i=1}^n \alpha_if_i$, where $\alpha_i \in B_1(X)$ and $f_i \in I$ for all $i=1,2,..,n$. For each $x \in X$, $\exists$ $k \in \{1,2,...,n\}$, such that $f_k(x) \neq 0$, otherwise it contradicts that $\mathbf{1}= \sum\limits_{i=1}^n \alpha_if_i$. We consider the map $g(x) = \sum\limits_{i=1}^n f_i^2(x)$, $\forall$ $x \in X$. Clearly, $g \in I\subseteq C(x)$ and $g(x) \neq 0$, $\forall$ $x \in X$. So $g$ is a unit in $I$. i.e., $I = C(X)$ - a contradiction. Hence, $I^*$ is a proper ideal of $B_1(X)$.
\end{proof}
\begin{theorem}
If $M$ is a hyperreal maximal ideal of $C(X)$ then $m(M^*)$ is a hyperreal maximal ideal of $B_1(X)$.   
\end{theorem}
\begin{proof}
If $m(M^*)$ is a real maximal ideal of $B_1(X)$ then by Theorem~\ref{RealMax3}, $m(M^*) \cap C(X)$ is a real maximal ideal of $C(X)$. Since $M \subseteq m(M^*) \cap C(X)$ and $M$ is maximal it follows that $M = m(M^*) \cap C(X)$ - a contradiction to the fact that $M$ is hyperreal. 
\end{proof}
\begin{theorem}\label{HyperReal}
The function $\zeta :  \mathcal{H}(C(X)) \ra \mathcal{H}(B_1(X))$ given by $\zeta (M) = m(M^*)$ is an injective function.
\end{theorem}
\begin{proof}
Let $M, N \in \mathcal{H}(C(X))$ be such that $m(M^*)= m(N^*)$. Then by maximality of $M$ and $N$ it follows that $M = m(M^*) \cap C(X) = m(N^*) \cap C(X) = N$.  
\end{proof}
\begin{corollary}\label{CorCard}
$|\MMM(C(X))| \leq |\MMM(B_1(X))|$.
\end{corollary}
\begin{proof}
This is immediate from Theorem~\ref{bij} and Theorem~\ref{HyperReal}.
\end{proof}

\section{Characterization of Real compact spaces }
\noindent From the discussion of the last section it follows that there is a one-one correspondence between the collections $\RRR\MMM(C(X))$ and $\RRR\MMM(B_1(X))$ given by $M \mapsto M_B$. It is well known \cite{GJ} that a Tychonoff space $X$ is real compact if and only if every real maximal ideal of $C(X)$ is fixed. Utilizing the one-one correspondence as mentioned above, we get a characterization of real compact spaces via real maximal ideals of $B_1(X)$.
\begin{theorem} \label{fixed}
A Tychonoff space $X$ is real compact if and only if every real maximal ideal of $B_1(X)$ is fixed.
\end{theorem}
\begin{proof}
Let $X$ be a real compact space and $\widehat{M} \in \RRR\MMM(B_1(X))$. By Theorem~\ref{bij}, there exists $M \in \RRR\MMM(C(X))$ such that $\widehat{M} = M_B$. Since $X$ is real compact, $M$ is fixed; i.e., $M = M_p$, for some $p \in X$. Hence, $\widehat{M} = M_B = {(M_p)}_B = \widehat{M}_p$ (by Example~\ref{FixedMax1}). \\
Conversely, let $M$ be any real maximal ideal of $C(X)$. Then $M_B \in \RRR\MMM(B_1(X))$ and so, $M_B$ is fixed. Therefore, $M$ ($ \subseteq M_B$) is a fixed ideal. Hence, $X$ is real compact.
\end{proof}
\noindent In Theorem 3.9 of our paper \cite{AA2}, we have proved a result for perfectly normal $T_1$-spaces though the same proof holds true for a larger class of spaces. In the following lemma, we state the result for the bigger class of spaces without proof. 
\begin{lemma}\label{lemma1}
If $X$ is a $T_4$-space in which every point is a $G_\delta$ point then the following statements are equivalent:
	\begin{enumerate}
	\item $X$ is finite.
	\item Every maximal ideal in $B_1(X)$ is fixed.
	\item Every ideal in $B_1(X)$ is fixed.
	
\end{enumerate}
\end{lemma}

\begin{theorem}
Let $X$ be a $T_4$ real compact space in which every point is a $G_\delta$-point. Then $B_1(X) = {B_1}^*(X)$ if and only if $X$ is finite.
\end{theorem}
\begin{proof}
If $X$ is finite then certainly, $B_1(X) = B_1^{*}(X)$. \\
Conversely, let $\widehat{M}$ be any maximal ideal of $B_1(X) = B_1^{*}(X)$. By Theorem 4.21 of \cite{AA2}, $\widehat{M}$ is a real maximal ideal. Since $X$ is real compact, by Theorem  \ref{fixed} $\widehat{M}$ is fixed. Finally, using Lemma~\ref{lemma1} we can conclude that $X$ is finite.
\end{proof}
\begin{corollary}
Let $X$ be a perfectly normal $T_1$ real compact space. $B_1(X) = B_1^{*}(X)$ if and only if $X$ is finite.
\end{corollary}

\end{document}